# Coreo-Graph: un progetto di teoria dei grafi e danza

Teresa Calogera*[1,2] e Tiziano Penati†[1,3]

[1]Dipartimento di Matematica “F. Enriques”, Università di Milano

[2]Istituto di Istruzione e Formazione Professionale "Immaginazione e Lavoro" - Sesto S. Giovanni (MI)

[3]GNFM (Gr. Naz. di Fisica Matematica) - Indam (Ist. Naz. di Alta Matematica “F. Severi”)

**Sunto /** In questo contributo rivisitiamo un laboratorio didattico di teoria dei grafi descritto in [Gaio & al. (2020)], basato sul problema dei ponti di Königsberg di Eulero, estendendone alcuni aspetti concettuali ed alcune attività laboratoriali; inoltre, proponiamo un'attività motoria ispirata alla danza, allo scopo di esplorare i cammini Euleriani e le loro proprietà. In questa nuova veste, il progetto è stato proposto e sperimentato nel triennio 2022-2024 come attività di rafforzamento delle competenze curriculari di *problem-solving*. Allo stesso tempo, esso è stato un'occasione di riflessione sul ruolo delle attività *embodied* nel facilitare la compresione ed assimilazione di concetti astratti della Matematica.

## 1 Introduzione

Il progetto qui descritto nasce con l'intento di offrire un laboratorio di Teoria dei Grafi che sia di potenziamento delle competenze di *problem-solving* nella Scuola Secondaria di Primo Grado (SSPG nel seguito) e nel biennio del Liceo Matematico (LM nel seguito), riprendendo la sperimentazione di *outdoor learning* descritta in [Gaio & al. (2020)] e valutandone possibili integrazioni e modifiche.

Una delle strategie consolidate nel proporre la Teoria dei Grafi nei vari ordini scolastici, come attività di potenziamento del pensiero critico e di *problem-solving*, è quella di presentarla come strumento per la risoluzione del celebre problema dei ponti di Königsberg di Eulero [Robinson (2006),Ferrarello & Mammana (2017), Santoso (2018), Gaio & al. (2020), Gonzales & al. (2021)]. Il problema (**Figura 1**, pannello di sinistra) consiste nel chiedersi se esista o meno un percorso che, partendo da una qualunque delle quattro regioni delimitate dal fiume Pregel, permetta di attraversare in sequenza tutti i sette ponti che collegano tali regioni, senza attraversare due volte lo stesso ponte. La non esistenza di un simile percorso si ottiene applicando l'enunciato del Teorema 1 (Teorem di Eulero) al grafo corripondente al problema (**Figura 1**, pannello di destra).

*teresa.calogera@gmail.com
†*corr. author:* tiziano.penati@unimi.it

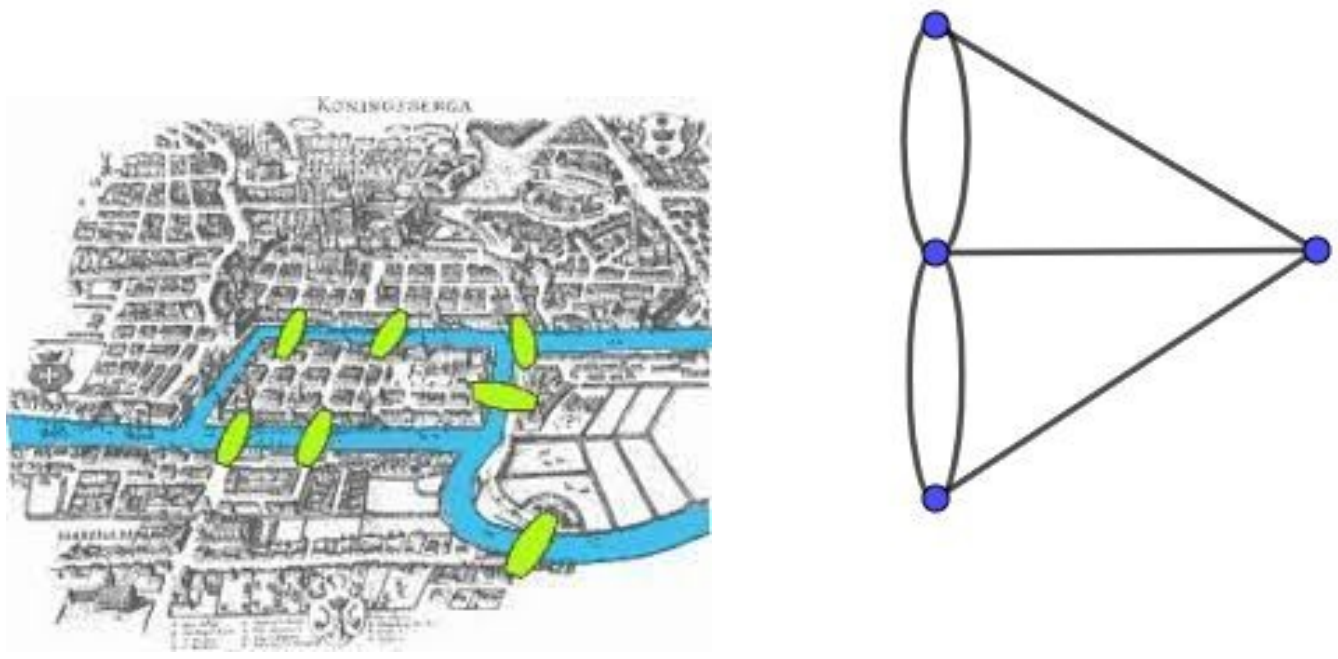

**Figura 1.** Illustrazione del problema dei ponti di Eulero (fonte `https://it.wikipedia.org/Problema_dei_ponti_di_Königsberg`) e del corrispondente grafo.

La proposta didattica contenuta in [Gaio & al. (2020)] è di riproporre il tema del problema dei ponti di Königsberg in un contesto motorio di orientamento all'interno di una mappa cittadina: l'idea dell'attività è quella di cercare, se esiste, una soluzione al *problema della guida turistica*, che deve accompagnare un gruppo di turisti lungo un itinerario, facendo in modo che ciascuna strada (che compone l'itinerario) venga percorsa una e una sola volta. Successivamente, la stessa tipologia di problema viene riproposta col linguaggio della teoria dei grafi, per mettere in luce la struttura geometrica del problema, il ruolo del grado di un vertice e illustrare il Teorema di Eulero sull'esistenza dei cammini Euleriani. La scelta metodologica è quella di integrare le attività laboratoriali di problem-solving con esperienze di tipo motorio, che rafforzino il riconoscimento della struttura di un problema matematico assegnato e delle sue proprietà, ancorandole così a un'esperienza concreta che ne faciliti il successivo studio astratto.

Il nostro progetto riprende questa scelta fatta dagli autori in [Gaio & al. (2020)]: abbiamo infatti incluso un'attività motoria ispirata alla danza, come fase di passaggio dalle attività iniziali sul problema dei ponti di Königsberg a quelle conclusive sulla teoria dei grafi. In tal modo, ci siamo riallacciati al quadro teorico e bibliografico dell'*embodied cognition*, che riconosce nelle esperienze motorio-sensoriali degli strumenti capaci di costruire riferimenti utili alla formazione e comprensione di concetti astratti [Nunez (2000), Nunez & al. (1999), Lakoff & Nunez (2000), Alibali & Nathan (2012)]. Nella nostra doppia sperimentazione abbiamo però potuto proporre tale esperienza solo all'interno del contesto della SSPG, mentre le classi del LM hanno dovuto svolgere le attività dei ponti e dei grafi in sequenza. Questa differenza nella presenza di un'attività motoria all'interno della sperimentazione, dovuta inizialmente a vincoli organizzativi, ci ha comunque spinto ad interrogarci sull'efficacia di un *laboratorio embodied* nella comprensione di nuovi concetti matematici, al netto delle evidenti differenze tra i due contesti con cui abbiamo avuto l'opportunità di interagire.

Un altro punto di contatto con la proposta descritta in [Gaio & al. (2020)] è stata la scelta di mantenere un percorso didattico che partisse dall'esplorazione di aspetti globali (cammini su mappe e grafi), per poi muovere verso gli aspetti più locali (il grado di un vertice) e infine alla sintesi dei due con la formulazione del Teorema di Eulero: questo sviluppo ricalca altri esempi in letteratura dove vengono proposte o descritte alcune sperimentazioni della teoria dei grafi nelle scuole [Ferrarello & Mammana (2017), Gonzales & al. (2021)].

A differenza di quanto descritto in [Gaio & al. (2020)], nella nostra progettazione e sperimentazione abbiamo quindi deciso di ampliare e modificare le attività laboratoriali in due direzioni:

- abbiamo aggiunto un costante focus, durante tutte le fasi del progetto, sulle proprietà di apertura e chiusura (circuiti) delle soluzioni, per rinforzare la riflessione sul ruolo del grado di un vertice all'interno dell'enunciato del Teorema di Eulero;
- abbiamo ripensato l'attività motoria per renderla più facilmente adattabile agli spazi disponibili nelle varie realtà scolastico-territoriali e per riportarla a una dimensione spaziale più confrontabile col corpo di chi la esegue; l'abbiamo resa in tal modo osservabile direttamente da parte dei docenti/esperti durante tutto il suo svolgimento.

In particolare, abbiamo esteso la prima attività sui ponti di Königsberg aggiungendo nuove mappe da esplorare che fossero differenti per tipologia di soluzioni, in modo da coprire tutte e tre le casistiche del Teorema di Eulero; le classi hanno così potuto lavorare fin dall'inizio sulle proprietà delle soluzioni e sul ruolo della regione di partenza. Nel progettare tale estensione ci siamo ispirati all'attività disponibile sul sito Mathigon [1], mantenendo la scelta di iniziare il laboratorio con la sfida più impegnativa, costituita proprio dal problema dei ponti di Königsberg; in tal modo le classi si sono subito confrontate col "dilemma" della possibile non-esistenza della soluzione al problema assegnato, a causa dall'apparente fallimento dei ripetuti tentativi effettuati. Questo ci è sembrato infatti il modo migliore per stimolare fin dall'inizio la creazione di quelle strategie e rappresentazioni personali del problema che sarebbero risultate utili allo svolgimento e alla comprensione delle attività successive [Kapur (2010), Kapur (2011), Jonsson & al. (2014)].

[1] `https://it.mathigon.org/course/graph-theory/bridges`

Per offrire invece un'attività motoria che fosse allo stesso tempo osservabile in ambito scolastico dai docenti/esperti e replicabile in qualunque contesto scolastico, ci siamo ispirati alla danza. Abbiamo sostituito il problema della guida turistica con quello che abbiamo chiamato *il problema delle coreografie Euleriane*, da cui il nome *Coreo-Graph*: l'obiettivo è quello di eseguire una sequenza di passi tra le posizioni assegnate (coreografia), che effettui ogni passo indicato esattamente una e una sola volta (si veda Sez. 4.2). Nonostante la difficoltà aggiuntiva di dover *fondere spazio e tempo* sincronizzando uno schema motorio con la base musicale [Gerofsky (2013)], questo laboratorio di danza consente di creare un'esperienza facilmente richiamabile nella successiva fase teorica, perchè rafforza e concretizza l'idea di cammino Euleriano su un grafo come un processo dinamico. La possibile riduzione del livello di astrazione [Hazzan & Hadar (2005)] che si ottiene da questa attività è stata bilanciata dalla richiesta di produrre sempre una rappresentazione simbolica di tali coreografie (stringhe alfanumeriche in corrispondenza biunivoca con le coreografie), proponendo contemporaneamente una visione più astratta del concetto di cammino su un grafo (sequenza di vertici e archi), e offrendo uno strumento linguistico efficace nell'esplorazione delle soluzioni al problema.

Il resto del lavoro è organizzato come segue. Nella Sezione 2 si richiamano sinteticamente i temi teorici portanti del progetto, ovvero il *problem-solving* e *l'embodiment*. Nella Sezione 3 vengono descritti gli obiettivi formativi e la metodologia adottata nello svolgimento del laboratorio. La Sezione 4 è dedicata alla descrizione del progetto nelle sue varie fasi e alla presentazione dell'attività *Coreo-Graph*. La Sezione 5 raccoglie i risultati della nostra sperimentazione, mentre la Sezione 6 è dedicata alle conclusioni e alle prospettive di sviluppo del progetto che abbiamo qui descritto.

## 2 Problem solving e embodiment

Il progetto, nel suo complesso, è un vero e proprio laboratorio di matematica, secondo la definizione che ne viene data in [AA. & VV. (2003)] : esso è infatti *«un insieme strutturato di attività volte alla* *costruzione di significati degli oggetti matematici»*, costruzione che *«è strettamente legata, da una parte, all'uso degli strumenti utilizzati nelle varie attività, dall'altra, alle interazioni tra le persone che si sviluppano durante l'esercizio di tali attività.»* Nell'illustrare il quadro teorico di riferimento, si delineano i diversi aspetti che emergono da tale definizione, esplicitando dunque quali sono i significati che vengono costruiti, gli strumenti attraverso cui essi vengono indagati e indicando la modalità di interazione tra esperti e studenti e quella da incoraggiare all'interno del gruppo classe.

## 2.1 Il problem solving

La comprensione degli elementi fondamentali della teoria dei grafi e il teorema di Eulero non costituiscono, dal punto di vista didattico, il fine del laboratorio proposto, quanto più un mezzo per mostrare e far sperimentare un'attività di problem solving. Il punto di partenza è la convinzione che, come afferma Alan H. Schoenfeld (si veda [Schoenfeld (1992)]) *«Mathematics instruction should provide students with a sense of the discipline - a sense of its scope, power, uses, and history. It should give them a sense of what mathematics is and how it is done, at a level appropriate for the students to experience and understand.»*

Le attività proposte offrono infatti una visione della matematica tutt'altro che procedurale. Lo studente, infatti, per portare a termine ciascuno dei compiti che gli vengono assegnati, non deve riprodurre schemi o algoritmi spiegati dall'insegnante, ma, piuttosto, è invitato a confrontarsi con un problema nuovo per la risoluzione del quale egli ha bisogno dell'aiuto di qualcuno di appena più esperto. Perciò, l'insegnante non si pone come unico detentore del sapere che fornisce istruzioni da seguire, quanto più come autentico soggetto che indaga e affronta il problema insieme agli studenti, con l'obiettivo di aiutarli a *«rendere il pensiero visibile»* [Collins & al. (1989)]. In questo modo, avviene un processo in cui l'attività cognitiva degli studenti e quella metacognitiva incoraggiata dall'insegnante si integrano continuamente e permettono non solo la comprensione del particolare, ma anche la scoperta - o l'approfondimento - di un modo di agire tipicamente matematico. Inoltre, come in [Gaio & al. (2020)], nella fase centrale del progetto, si è scelto di proporre un problema del tutto analogo a quello dei ponti, ma con un radicale cambio di contesto. Ciò permette di mostrare che un cambio di prospettiva può essere d'aiuto alla comprensione e alla risoluzione di un problema, cosa ben nota ai matematici ma che difficilmente si riesce a mostrare sui banchi di scuola.

In sintesi, le attività proposte rappresentano un esemplare approccio al problem solving, poichè permettono di ripercorrere e conquistare i passi necessari alla comprensione e risoluzione di un problema posto (e risolto) da uno dei più grandi matematici della storia, Eulero.

## 2.2 L'embodied cognition

Lo strumento scelto per avvicinarsi alla risoluzione del problema è un'attività *embodied* che utilizza la danza, con la convinzione che l'azione del corpo possa sostenere la formulazione del pensiero e facilitare così il passaggio dal problema reale alla teoria astratta che si sviluppa per risolverlo [Lakoff & Nunez (2000)].

In particolare, nel caso in esame, la delicatezza della traduzione dal problema reale al modello astratto sta tutta nella scelta dell'identificazione delle regioni con i vertici e dei ponti con gli archi del grafo. Per far sì che gli studenti colgano tale aspetto, è innanzitutto fondamentale prestare attenzione al lessico utilizzato durante la prima attività: il problema di Königsberg riguarda l'attraversamento dei ponti e non delle regioni, che sono viste come punti di arrivo e di partenza dei ponti stessi. L'introduzione di una traduzione intermedia tra tali mappe e i grafi, attraverso l'utilizzo degli schemi coreografici, permette di ridurre i possibili fraintendimenti non tanto perchè richiede una riflessione più approfondita, quanto perchè il pensiero diventa, grazie a questa attività, azione del corpo: i ponti da attraversare tra una regione e l'altra diventano passi da svolgere tra due posizioni nello spazio.

La richiesta di percorrere fisicamente le coreografie euleriane può rendere inoltre più evidenti le proprietà delle stesse (del tutto analoghe a quelle delle soluzioni del problema dei ponti), anche agli occhi di chi non è riuscito ad acquisire familiarità con la notazione simbolica introdotta nell'attività delle mappe e dei ponti. Infatti, l'osservazione di tali proprietà che prima poteva essere frutto di un ragionamento sulla geometria dei cammini o sulle stringhe di simboli, diventa qui semplice presa di coscienza di un'azione fisica.

Questo è un vantaggio dell'aver disegnato un'attività motoria che si svolge all'interno del mesospazio, a differenza di quella proposta in [Gaio & al. (2020)] che chiede di muoversi in un macrospazio. Il corpo si ritrova così immerso in uno spazio confrontabile con le proprie dimensioni e lo esplora riuscendo a cogliere la struttura delle coreografie eseguite con un ridotto livello di concettualizzazione rispetto al dover compiere un percorso nello spazio urbano. Inoltre, con tale scelta i gruppi hanno sempre una visione duplice del problema; infatti, mentre chi balla sullo schema ne percepisce la struttura dall'interno, i compagni che lo osservano dall'esterno mantengono una prospettiva globale. Questo doppio punto di vista aiuta a rendere più evidente l'esistenza di eventuali proprietà locali e globali del grafo che si riflettono nella struttura delle coreografie euleriane.

Infine, un'osservazione a proposito del fatto che l'attività si svolge all'interno del mesospazio, a differenza di quella proposta in [Gaio & al. (2020)] che chiedeva agli alunni di muoversi in un macrospazio. Tale scelta comporta che gli studenti abbiano sempre una visione duplice del problema; infatti, mentre uno balla sullo schema e ne percepisce la struttura dall'interno, i compagni che lo osservano dall'esterno mantengono una prospettiva globale. Questo doppio punto di vista può rendere più evidente l'esistenza di eventuali simmetrie all'interno del grafo che si riflettono nella struttura delle coreografie euleriane.

# 3 Obiettivi, metodi e strumenti

## 3.1 Contesto di riferimento e durata dell'attività

Il laboratorio qui descritto è stato rivolto principalmente alla Scuola Secondaria di primo grado e, in una formula opportunamente adattata, al biennio delle Scuole Secondarie di secondo grado. Nel primo caso, il progetto ha coinvolto 9 classi seconde della SSPG durante il triennio 2022-2024 di sperimentazione, per un totale di oltre 150 partecipanti, ed è stato diviso in tre incontri in orario curricolare, della durata rispettivamente di due ore (Fase 1, problema dei ponti), un'ora (Fase 2, *Coreo-Graph*, ovvero ricerca di coreografie Euleriane) e due ore (Fase 3, teoria dei grafi, Teorema di Eulero e sue applicazioni). In questo contesto, abbiamo cercato di mantenere i tre incontri molto ravvicinati, per non disperdere la comprensione e l'attenzione acquisite durante ciascuna fase. Nel secondo caso, il progetto, privato della parte motoria inclusa nella seconda fase, ha coinvolto nello stesso triennio 5 classi del LM, per un totale di quasi 150 partecipanti, ed è stato svolto presso il Dipartimento di Matematica "F.Enriques" di Milano nell'arco di mezza giornata.

## 3.2 Metodologie didattiche e strumenti

Il progetto rispecchia la suddivisione già illustrata in [Gaio & al. (2020)] in tre fasi, che rappresentano tre successivi momenti di lavoro con ciascuna classe: in tale suddivisione, il problema dei ponti di Königsberg proposto inizialmente, e lasciato insoluto, viene ripreso solo alla fine dell'ultima fase e qui, opportunamente tradotto nel linguaggio della teoria dei grafi, viene risolto. In questa rivisitazione del progetto, abbiamo deciso di non limitarci alla questione della sola esistenza delle soluzioni, ma di estendere l'analisi delle loro proprietà, focalizzandoci sul ruolo svolto dai vertici di partenza. Infatti, in tutte e tre le fasi del progetto, abbiamo sempre proposto le stesse due domande:

1. esistono delle soluzioni al problema?
2. esistono delle soluzioni partendo da qualunque regione/posizione/vertice?

Questa estensione ha avuto un duplice scopo: anzitutto riconoscere la stessa struttura del problema nelle tre differenti situazioni, al fine di favorire la formulazione e comprensione dell'enunciato del Teorema di Eulero. Allo stesso tempo abbiamo voluto mostrare una metodologia che caratterizza l'agire matematico: indagare le proprietà delle soluzioni talvolta può offrire una visione più ampia e chiara circa la questione dell'esistenza delle stesse.

L'attività della danza, che rappresenta il vero ingrediente di originalità di questo progetto, ha avuto il vantaggio di poter essere facilmente adattata ai diversi contesti in cui è stata sperimentata (palestre, auditorium, giardini), e ha permesso non solo di osservare le classi al lavoro anche in questa fase, ma di stimolarle proponendo variazioni al problema assegnato (aggiunta/rimozione di un passo o composizione di due distinte coreografie in sequenza).

Il punto di arrivo del progetto è l'enunciato del Teorema dei cammini Euleriani (su grafi connessi): l'utilizzo di tale strumento per rileggere e risolvere tutte le situazioni precedenti (ponti e danza) aiuta a evidenziare la differenza tra risolvere un problema specifico per tentativi e cercare una strategia generale, che sia poi di facile applicazione e che accomuni tanti diversi problemi con la stessa struttura.

Nelle varie fasi sono stati proposti alle classi diversi strumenti e modalità per lo studio dei cammini Euleriani sui grafi connessi. Abbiamo proposto delle attività con carta e penna, facili da riproporre in ogni contesto scolastico; laddove possibile, abbiamo permesso alle classi di lavorare sui tablet degli Istituti (laboratori informatici itineranti), per poter accedere alle attività del sito Mathigon in modalità digitalmente sicura e controllata. Alcune fasi di restituzione post-attività sono state svolte con un Padlet condiviso, in modo da riprendere assieme in aula le considerazioni fatte singolarmente dagli studenti a casa. L'attività motoria delle coreografie Euleriane è stata invece introdotta da video esplicativi, per anticipare in aula la spiegazione del tipo di attività che avremmo chiesto loro di svolgere successivamente.

### 3.3 Obiettivi formativi

Come già osservato da Mammana e Milone (2010) nella loro sperimentazione con analoghe fasce d'età, l'assenza di prerequisiti specifici nell'ambito del curriculum di matematica consente alla maggior parte dei partecipanti di potersi mettere in gioco, servendosi principalmente della propria curiosità e delle proprie capacità logico-intuitive.

Stimolando questi due ambiti, il progetto vuole operare allo stesso tempo su due differenti aspetti della pratica laboratoriale. Da una parte vuole cercare di superare lo scetticismo, mediamente diffuso tra gli studenti della SSPG, verso l'utilità e la necessità dei contenuti della disciplina (a differenza di quanto solitamente accade con molti altri contenuti dell'area STEM), lavorando in maniera guidata sul processo di creazione degli strumenti matematici necessari alla risoluzione del problema assegnato ("Ha rafforzato un atteggiamento positivo rispetto alla matematica attraverso esperienze significative e ha capito come gli strumenti matematici appresi siano utili in molte situazioni per operare nella reltà" - Traguardi per lo sviluppo delle competenze, Zanichelli); dall'altra, esso vuole accompagnare le classi verso un'azione creativa che è difficile da includere sistematicamente nella didattica frontale della matematica, dove la teoria di norma precede la sua applicazione. In questo caso la pratica laboratoriale vuole accompagnare alla formulazione dell'enunciato del Teorema di Eulero, ovvero lo strumento teorico necessario per risolvere e reinterpretare a posteriori i problemi già incontrati, così come nuovi problemi in cui si riconosce la stessa struttura. In relazione a tali finalità, possiamo elencare questi obiettivi specifici:

- risolvere un problema concreto tramite un'opportuna traduzione nel linguaggio matematico (modellizzazione matematica e problem solving);
- passare da un approccio sperimentale ai problemi matematici alla successiva formulazione di una congettura sotto forma di enunciato (azione creativa);
- riflettere sull'esistenza e non esistenza delle soluzioni di un problema e imparare a operare sul problema inverso (ovvero come modificare un grafo/una mappa in modo da avere soluzioni con determinate proprietà) come metodologia utile alla risoluzione del problema diretto.

# 4 Descrizione del progetto

Descriviamo ora le tre fasi del progetto. Come già messo in evidenza nel paragrafo precedente, tutte le fasi sono accommunate dalle stesse domande guida, opportunamente adattate al problema specifico della fase in svolgimento:

1. esiste una soluzione al problema assegnato?

2. riesco a trovare una soluzione partendo da una qualsiasi regione/posizione/vertice?

L'obiettivo è mettere in luce la struttura che accomuna tutti i problemi proposti, per poter affrontare in maniera più efficace l'enunciato del Teorema di Eulero, come punto di arrivo del percorso di *problem-solving*.

## 4.1 FASE 1: I ponti di Eulero e le sue varianti

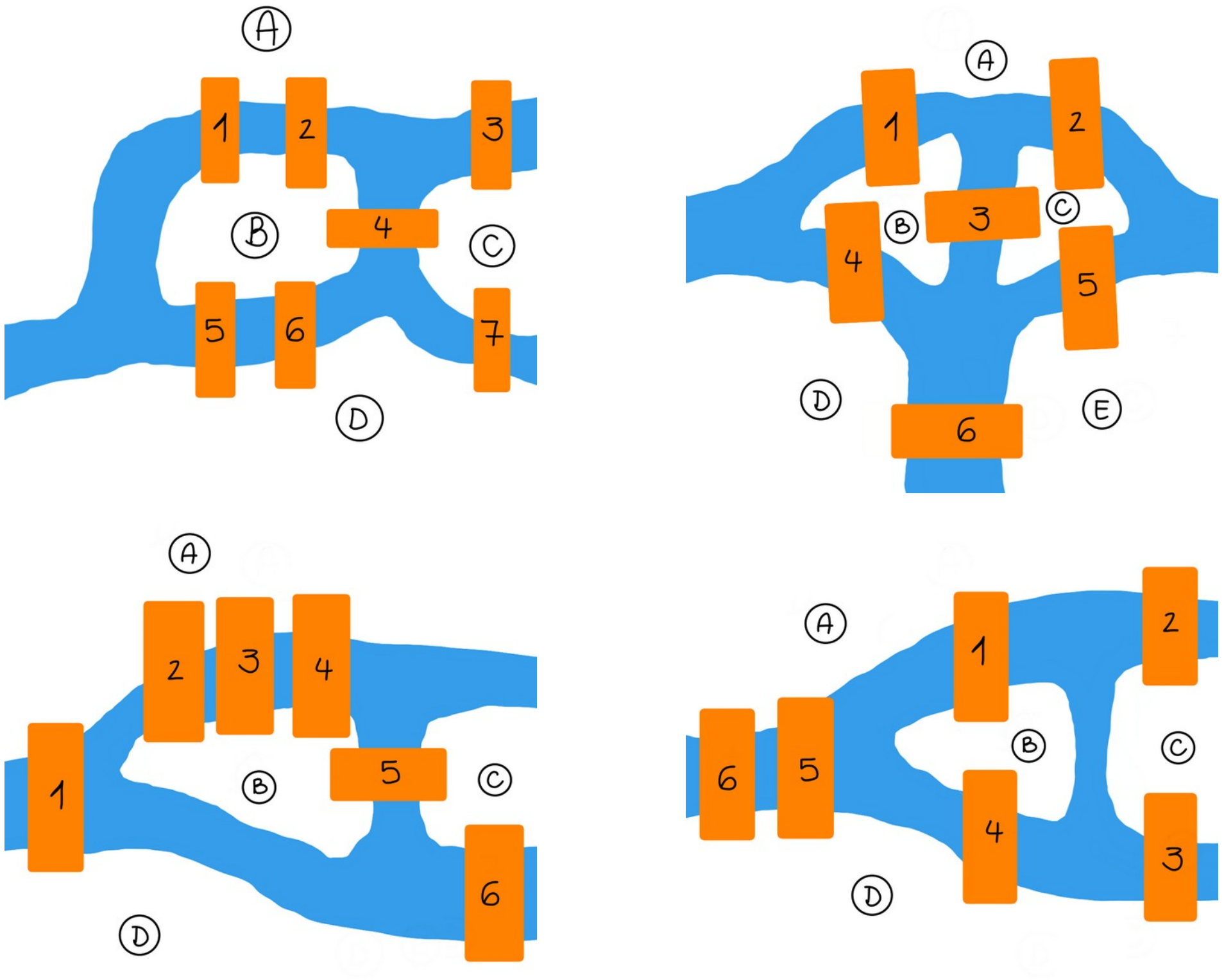

**Figura 2.** Mappe usate per la prima fase. In alto, riproduzione delle prime due mappe presenti dal sito Mathigon (la prima corrispondende al problema di Eulero originale). In basso due mappe usate come esempi dell'esistenza di soluzioni a partire da ogni regione, studiate attraverso la rappresentazione simbolica.

La prima fase si focalizza sulla risoluzione del problema dei ponti, nella versione originaria del problema di Eulero e in alcune sue varianti. Abbiamo chiamato *percorso Euleriano* (o cammino Euleriano sulla mappa) una soluzione a tale problema, ovvero un percorso che, partendo da una regione, attraversa una e una sola volta tutti i ponti della mappa. Per questa fase ci siamo inizialmente appog-giati e ispirati alle attività proposte sul sito Mathigon: il sito permette di confrontarsi col problema dell'esistenza e non esistenza delle soluzioni con uno strumento digitale che facilita la ricerca delle soluzioni per tentativi (offrendo la possibiltà di ripulire la mappa dopo ogni tentativo). In aggiunta alle mappe già proposte da Mathigon, abbiamo inserito nuove mappe che includono il caso non contemplato, ovvero quello di mappe che ammettono soluzioni a partire da qualunque regione (il caso dei circuiti). La **Figura 2** mostra alcune possibili mappe, nell'ordine con cui sono state sempre assegnate.

La scelta di assegnare inizialmente la mappa più complessa, ovvero quella del problema originale di Eulero, permette di far partire il laboratorio con una narrazione ancorata al contesto storico del problema; allo stesso tempo, immerge fin da subito le classi nella difficoltà del problema e nella fatica necessaria per elaborare una buona strategia risolutiva.

Affrontate le prime due mappe (in alto, **Figura 2**), abbiamo introdotto la *rappresentazione simbolica REGIONI (lettere)-PONTI (numeri interi)* di un percorso sulla mappa. Ad esempio, nella mappa in basso a sinistra in **Figura 2**, si possono numerare i ponti da sinistra verso destra: una possibile soluzione che inizia e termina nella regione $A$ si scriverà quindi come

$$A1D6C5B4A3B2A\,. \tag{1}$$

In questo modo abbiamo offerto uno strumento simbolico che aiuta a tener traccia delle soluzioni trovate e consente allo stesso tempo di indagarne alcune proprietà in modo più astratto: per esempio, osservando che ogni permutazione delle terne regione-ponte-regione fornisce nuove e distinte soluzioni[2].

### 4.1.1 Discussione e conclusioni

La breve discussione finale di questa prima fase pone l'accento sul confronto tra le tre situazioni presentate, che esemplificano i tre possibili casi del Teorema di Eulero. In questa fase vogliamo anzitutto far emergere le differenze tra le tre situazioni, per facilitare l'analogia con la fase successiva (problema della coreografia Euleriana) e per iniziare a discutere delle *strategie di problem solving*. Abbiamo sempre concluso assieme alle classi che:

[2] Sfruttando il fatto che le soluzioni sono "circuiti" (stessa regione di partenza e arrivo), si può spostare la terna regione-ponte-regione $A1D$ alla fine del percorso e ottenere (senza ripetere due volte la regione $A$) la soluzione

$$D6C5B4A3B2A1D\,. \tag{2}$$

che ha come nuova regione di partenza e arrivo la regione $D$.

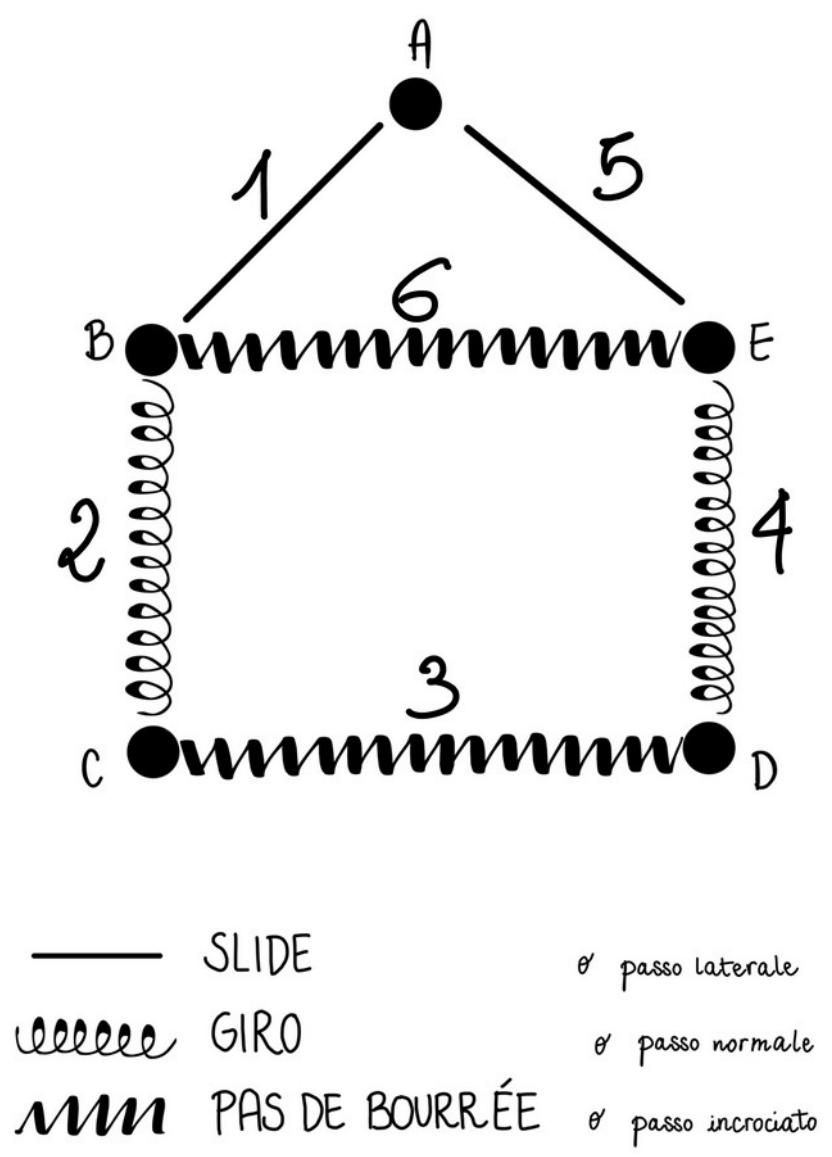

**Figura 3.** Schema coreografico a "casetta" con numerazione dei passi.

- la risolubilità di un problema è stata sempre ottenuta attraverso la scoperta di almeno una soluzione. Mostrare la non risolubilità è invece molto più difficile.
- Infatti, quando non riusciamo a trovare la soluzione, non siamo sicuri della sua esistenza o non esistenza, a meno di esplorare tutti i cammini possibili. Andare per tentativi è faticoso e poco efficiente come approccio al "problem solving", ma anche necessario per costruirsi una prima rappresentazione della situazione matematica da indagare.

## 4.2 FASE 2: il problema della *coreografia Euleriana*

Col problema della *coreografia Euleriana* vogliamo riproporre il problema della ricerca di cammini Euleriani su un grafo nel contesto della danza: abbiamo cercato di progettare un'attività capace di combinare un'analisi geometrica della situazione proposta con una disciplina artistico-motoria in grado di intercettare alcuni interessi o tendenze della fascia d'età con cui abbiamo lavorato.

La prima parte di questa attività si svolge in classe. Il punto di partenza è lo schema di **Figura 3**, di cui abbiamo precisato gli elementi costitutivi, anticipando di fatto la struttura del grafo (esattamente come era stato anticipato in [Gaio & al. (2020)] dall'itinerario sulla mappa):

- i vertici - ovvero le posizioni in cui si può trovare il ballerino;

- i passi - ovvero le transizioni che sono consentite;

- lo stile dei passi - ovvero il modo con cui si effettua ciascuna transizione.

Abbiamo introdotto per chiarezza una definizione di *coreografia* sullo schema assegnato, in analogia al concetto di percorso sulle mappe della prima fase: con *coreografia* si intende ogni sequenza di passi, che parta e termini nei vertici, compatibile con lo schema. Il ballerino può quindi spostarsi solo tra vertici collegati da passi (a prescindere dallo stile del passo). Allo stesso modo, ripensando al problema dei ponti (i ponti sono ciò che ci fa *passare* da una regione all'altra), diventa naturale definire una *coreografia Euleriana* come una sequenza di passi dello schema assegnato che consenta al ballerino di compiere tutti i passi dello schema, ma facendo ogni passo una e una sola volta. Servendoci anche in questo contesto della rappresentazione simbolica già introdotta, possiamo fornire quindi un esempio di coreografia Euleriana: $B6E5A1B2C3D4E$. Per rinforzare il concetto di *euleriani-tà*, durante la sperimentazione abbiamo mostrato dei video autoprodotti: alcuni in cui viene eseguita una coreografia euleriana e altri in cui viene commesso un errore (un passo non viene eseguito/viene eseguito due volte).

#### 4.2.1 L'attività *Coreo-Graph*

Dopo aver diviso gli studenti in gruppi (o crew), assegniamo a ogni gruppo uno schema di danza equivalente a un grafo, come quello in **Figura 3**. Abbiamo distribuito tre diverse tipologie di schemi (uno schema per ciascuna crew), di cui due risolubili (nel senso dell'esistenza di coreografie Euleriane) e uno non risolubile:

- uno schema risolubile con (due) vertici di grado dispari (GC1);

- uno schema risolubile con soli vertici di grado pari (GC2);

- uno schema non risolubile (GC3).

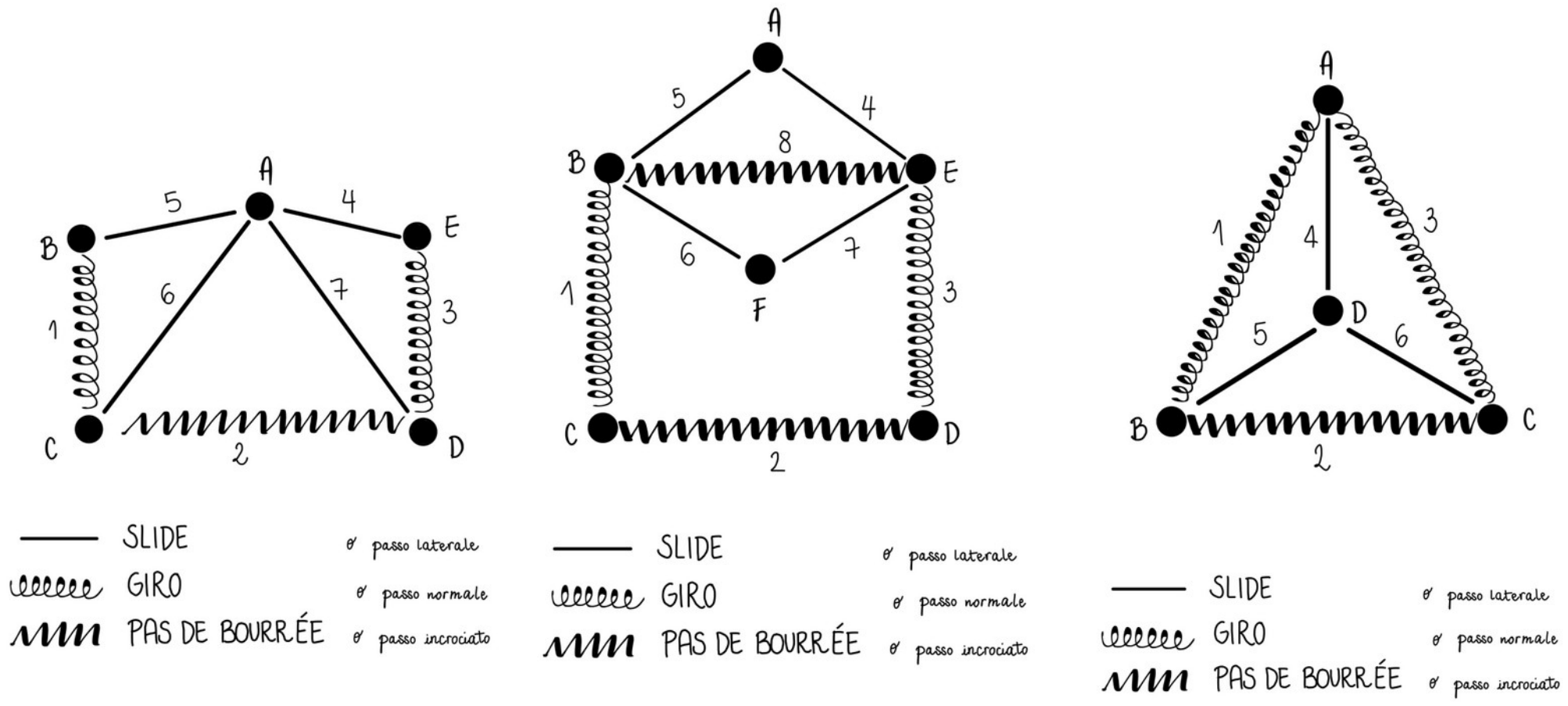

**Figura 4.** In ordine da sinistra a destra: CG1, CG2, CG3.

Abbiamo poi reso visibili gli schemi sul pavimento attraverso dei "bullets" adesivi colorati per le posizioni e dei nastri adeisivi per i passi (come in **Figura 5**).

Al fine di lavorare sulle stesse due domande guida della fase precedente, opportunamente declinate nella nuova attività, abbiamo chiesto alle crew di:

1. identificare possibili coreografie Euleriane e scriverle con la rappresentazione simbolica, prestando attenzione al ruolo svolto dalla scelta del vertice di partenza (posizione iniziale del ballerino). Dopo una prima fase di esplorazione, alla crew in possesso dello schema GC3 viene proposto di aggiungere un passo tra due posizioni a scelta, in maniera tale da assicurare l'esistenza di almeno una coreografia Euleriana.

2. eseguire (a turno) la coreografia senza base musicale, usando i riferimenti sul pavimento;

3. scegliere due elementi della crew disposti a memorizzare le coreografie sincronizzandosi con la base musicale, eventualmente modificando a proprio piacimento lo stile dei passi indicati nello

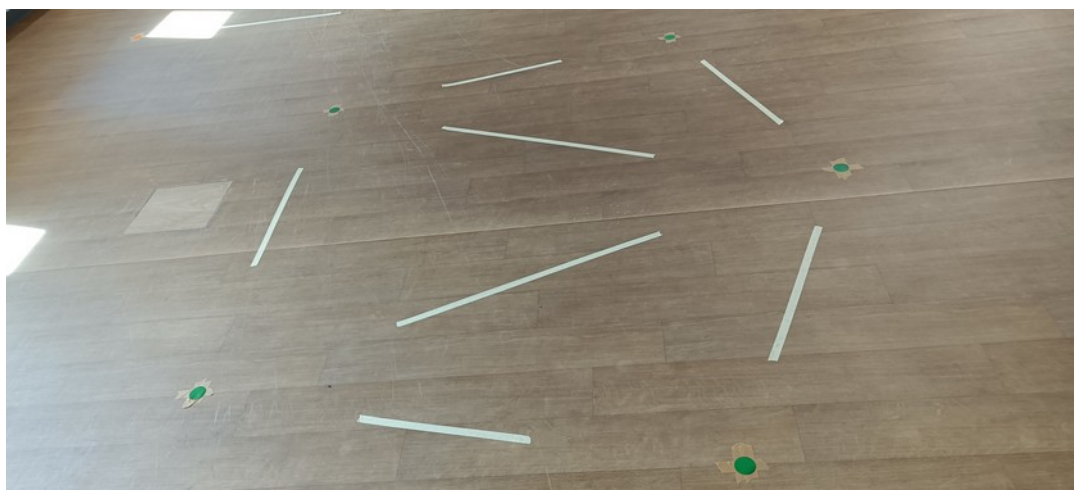

**Figura 5.** Schema CG1 riprodotto in palestra.

schema;

4. far eseguire ai/alle ballerini/e selezionati dalle crew le proprie coreografie Euleriane, nella breve esibizione finale.

#### 4.2.2 Discussione e conclusioni

A conclusione di questa attività, abbiamo proposto alle crew di “postare” le proprie risposte alle due domande guida su un Padlet condiviso e amministrato dagli esperti esterni, invitando ciascuna crew a osservare gli schemi e le soluzioni degli altri gruppi. Questa fase di restituzione su Padlet verrà poi ripresa e commentata in classe durante l’inizio della terza e ultima fase, prima di presentare la teoria dei grafi.

### 4.3 FASE 3: introduzione alla teoria dei grafi e Teorema di Eulero

L’inizio della terza fase riprende il problema della coreografia Euleriana sviluppato sul Padlet dalle varie crew, al fine di far emergere con chiarezza le diverse proprietà delle coreografie Euleriane, nei due casi in cui queste esistono:

1. quando esistono coreografie euleriane che partono da ogni vertice dello schema, la coreografia è sempre un circuito;

2. quando invece i possibili vertici di partenza per le coreografie sono solamente due, essi sono sempre il punto di partenza e di arrivo di una coreografia Euleriana “aperta”.

Prima di abbandonare la riflessione sull’attività motoria in favore della teoria dei grafi, abbiamo esteso agli schemi CG1 e CG2 il *problema inverso* assegnato nella seconda fase al solo schema GC3; ovvero, abbiamo chiesto a ogni crew di rimuovere un solo passo in maniera tale da modificare la tipologia delle proprie coreografie Euleriane (da chiuse ad aperte e viceversa). In tal modo, si incrementano le esperienze a cui riferirsi nella successiva riflessione sull’importanza della parità del grado dei vertici nella risoluzione del problema assegnato e nella tipologia delle sue soluzioni. La terza e ultima fase prosegue con due momenti distinti:

1. la ricerca di cammini Euleriani su grafi connessi e la classificazione dei grafi secondo le proprietà di tali cammini, fino alla formulazione del Teorema di Eulero;

| grafo | esistenza | dipendenza |
|---|---|---|
| G1 | No | nonsense |
| G2 | Si | Si |
| G3 | Si | Si |
| G4 | Si | No |
| G5 | No | nonsense |
| G6 | No | nonsense |

| tipo | grafi |
|---|---|
| I tipo | G4 |
| II tipo | G2, G3 |
| III tipo | G1, G5, G6 |

**Tabella 1.** Tabelle per la classificazione dei tre tipi di grafi.

2. l'utilizzo del Teorema di Eulero per rileggere il problema inverso delle coreografie euleriane e infine *tradurre e risolvere* il problema dei ponti di Königsberg e alcune sue possibili varianti.

### 4.3.1 Ricerca di cammini Euleriani su grafi connessi

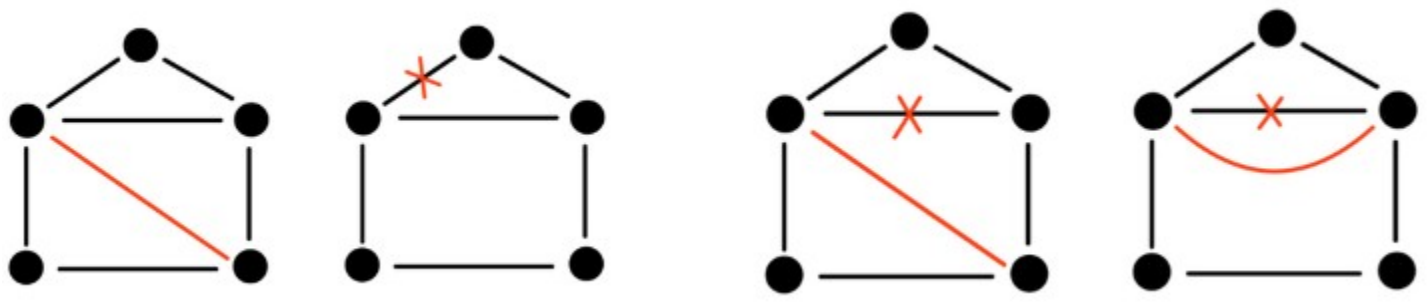

**Figura 6.** Trasformazioni. In ordine: aggiunta, rimozione, taglia e incolla

Dopo aver sottolineato come il numero di vertici, di passi e il modo con cui i passi collegano i vertici sia tutto ciò che caratterizza la risolubilità del problema della ricerca di una coreografia Euleriana negli schemi assegnati, introduciamo gli elementi di un grafo (vertici, archi, estremi di un arco) e ci soffermiamo sulla nozione di grado di un vertice. Abbiamo deciso di sottolineare maggiormente questa nozione, osservando come cambino i gradi dei cinque vertici nelle tre trasformazioni iniziali sul grafo (denominato in seguito G2) in **Figura 6**. Questo nella speranza di facilitare il riconoscimento del ruolo del grado dei vertici nei cammini Euleriani dei grafi assegnati successivamente.

L'allegato 1 mostra una possibile scheda usata per questa attività, guidata dalle solite due domande delle fasi precedenti:

1. esistono cammini Euleriani?

2. riusciamo a trovare cammini Euleriani partendo da qualunque vertice?

Alla fine di questa prima attività abbiamo deciso di formalizzare maggiormente la classificazione dei grafi rispetto alle due domande assegnate. Abbiamo raccolto i risultati in due tabelle: una prima tabella (doppia entrata) con le risposte alle due domande per ogni grafo, una seconda tabella con la classificazione dei problemi in funzione di tali risposte (primo, secondo e terzo tipo). In tal modo abbiamo già schematizzato le tre situazioni previste dall'enunciato del Teorema di Eulero.

La parte conclusiva di questo primo momento serve a formulare l'enunciato del Teorema di Eulero. Per provare a fare emergere il ruolo del grado nei tre differenti tipi di grafo, ci siamo anzitutto focalizzati sui grafi del secondo tipo: l'obiettivo è riconoscere il ruolo del grado dei due vertici "speciali" da cui iniziano e terminano tutti i cammini Euleriani. Per facilitare il riconoscimento del ruolo della (dis)parità del grado di tali vertici, è possibile far riflettere sul percorso che si ottiene quando si traccia un cammino Euleriano sul grafo G2, partendo da uno dei 2 vertici di grado dispari (come nel gioco del disegno della casetta senza staccare mai la matita dal foglio e senza rifare due volte lo stesso segmento), sviluppando un argomento di tipo *entry-exit* nei vertici che si incontrano lungo lo sviluppo del cammino (si veda anche [Ferrarello & Mammana (2017)]). Infatti, *dai vertici di partenza/arrivo si esce/entra una volta in più di quanto si esca/entri, mentre tutti gli altri vertici sono di passaggio, quindi si entra ed esce lo stesso numero di volte*. Per avere un vertice di partenza e uno di arrivo, questi devono necessariamente avere grado dispari, mentre tutti gli altri devono avere grado pari. Una volta riconosciuto il ruolo discriminante svolto dal grado dei due vertici in questa tipologia di grafo, si procede per confronto con gli altri due tipi di grafi, giungendo alla formulazione dell'enunciato del

**Teorema 1 (Teorema sui cammini/circuiti Euleriani in grafi connessi)** *In un grafo connesso $G$, di vertici $v_i \in V$:*

**I)** *se tutti i vertici $v_i$ hanno grado pari, allora esistono cammini Euleriani e non dipendono dalla scelta iniziale del vertice. Inoltre, il vertice di partenza e di arrivo coincidono (ovvero sono circuiti Euleriani).*

**II)** *se esistono solo due vertici $v_{1,2}$ di grado dispari, allora esistono solo cammini Euleriani che partono in $v_1$ e terminano in $v_2$ o viceversa (quindi non sono circuiti).*

**III)** *in tutti gli altri casi non esistono cammini Euleriani.*

Abbiamo concluso questa prima attività rileggendo il problema della coreografia Euleriana alla luce del Teorema di Eulero, soffermandoci sul cambio del grado dei vertici negli schemi a seguito di rimozioni (o aggiunte) di passi tra i vertici esistenti.

#### 4.3.2 Risoluzione del problema dei ponti di Königsberg

Il secondo momento inizia con la proposta di usare il Teorema di Eulero appena enunciato per risolvere il problema dei ponti di Königsberg; per collegare la mappa con il corrispondente grafo bisogna infatti identificare ogni regione con un punto (come se fosse un punto a caso da cui partire all'interno della stessa regione) e ogni ponte con un arco ( può risultare necessario esplicitare l'analogia tra gli archi che si possono attraversare una e una sola volta, e lo stesso ruolo svolto dai ponti). Effettuata tale trasformazione (per esempio con una animazione, nel nostro caso costruita con Geogebra), diventa naturale identificare la mappa di partenza con il grafo G1 della precedente attività (si veda l'allegato 1), concludendo la non esistenza delle soluzioni.

Dopo aver mostrato come esplorare le mappe della prima fase trasformandole in opportuni grafi, è possibile proporre delle sfide conclusive sulle mappe per valorizzare l'efficacia di un'opportuna modellizzazione matematica del problema e dello strumento teorico offerto dal Teorema di Eulero. Nel nostro caso, ci è sembrato interessante proporre la variante dell'abbattimento di un ponte e sua ricostruzione in modo da rendere risolvibile il problema («è possibile demolire e ricostruire il ponte sia rendendo il grafo del I tipo sia del II tipo?»), per analogia a quanto già visto nel contesto delle coreografie Euleriane. In aggiunta, abbiamo proposto di risolvere il problema della ricerca di percorsi Euleriani su una opportuna schematizzazione della mappa della città di Leiden (si veda, ad esempio, ), per rinforzare il contrasto tra l'approccio iniziale e la semplicità che ci viene offerta dal Teorema nel mostrare la non esistenza delle soluzioni (la schematizzazione scelta presentava almeno quattro regioni con "grado dispari", ovvero con un numero dispari di ponti).

## 5 Risultati della sperimentazione

Raccogliamo in questa sezione i principali risultati della sperimentazione del progetto. Abbiamo deciso di raggruppare i risultati per tematiche; dove ci è parso opportuno, abbiamo messo in luce le differenze tra la sperimentazione fatta nella SSPG e quella fatta nel biennio del LM.

## 5.1 Mappe, ponti e le strategie di problem solving

1. La scelta di iniziare il progetto con il problema dei ponti di Königsberg, senza aver fornito gli strumenti adeguati alla sua completa risoluzione, ha permesso a ciascuno studente di sviluppare la propria personale rappresentazione geometrica del problema e le proprie strategie risolutive, a prescindere dall'apparente insuccesso nell'utilizzo delle stesse [Kapur (2010), Kapur (2011), Jonsson & al. (2014)]. Ne sono stati la conferma da una parte la rapidità mostrata nel risolvere le successive attività sulle mappe, dall'altra i diversi commenti su quali elementi del problema sembrassero costituire l'ostruzione all'esistenza di un percorso Euleriano. A tal proposito, dalle classi sono sempre emerse considerazioni stimolanti sulla natura del problema: per alcuni l'ostruzione sembrava esser rappresentata *dal numero dispari di ponti* nella mappa iniziale di **Figura 2**, per altri *la mancanza di un ponte* in una specifica posizione della mappa o la *presenza di un ponte in eccesso* (entrambe manifestazioni di un "errato" grado della regione).

2. Di fronte al problema dell'esistenza di un percorso Euleriano nella mappa di Königsberg, tutte le classi si sono sempre sostanzialmente divise in due gruppi: da una parte chi era disposto a scommettere sull'esistenza di un percorso euleriano, che risultava però inaccessibile con la strategia usata, dall'altra chi era ragionevolmente convinto della sua non esistenza, a fronte dei numerosi tentativi effettuati. Nel caso del LM, questa seconda posizione era spesso sostenuta da argomenti (non sempre corretti) relativi alla particolare distribuzione dei ponti tra le regioni assegnate, o addirittura facendo riferimento alla disparità dei ponti costruiti in una specifica regione di partenza, di fatto anticipando una congettura sul ruolo del grado di un vertice nel corrispondente grafo. Viceversa, all'interno delle classi della SSPG abbiamo riscontrato che *il credere nell'esistenza di una soluzione*, più che una congettura supportata dalla sperimentazione, pareva esser frutto del contratto didattico secondo cui un docente di norma non propone un problema impossibile. Tuttavia, al termine dell'attività sulle tre mappe, una parte significativa degli studenti si è spesso ritrovata ad aver cambiato opinione sull'esistenza di un percorso Euleriano nel problema dei ponti di Königsberg, a favore della non esistenza.

3. La simmetria della seconda mappa in **Figura 2** ha ridotto la complessità del problema, consentendo di dedurre il ruolo svolto dalle regioni centrali con una esplorazione parziale di tutti i percorsi ammissibili: infatti, ogni percorso costruito a partire da una delle due regioni di destra si "riflette" in un corrispondente percorso che parte da una delle due regioni di sinistra, e viceversa. Questa strategia è stata adottata tanto nella SSPG quanto nel LM. Nella SSPG ha però rischiato di spostare l'attenzione su un aspetto secondario del problema: la simmetria è diventata in alcuni casi la proprietà discriminante per l'esistenza di percorsi Euleriani con inizio e fine dalle "isole centrali".

## 5.2 Aspetti linguistici

1. Nelle prime battute dell'attività sulle mappe abbiamo dovuto sempre sottolineare il fatto che i cammini dovessero *partire da una regione* e non da un ponte e che questo ultimo dovesse essere *attraversato interamente* da una regione all'altra (questo per evitare che venissero considerate come soluzioni quelle costruite attraversando solo parzialmente un ponte). In tal modo abbiamo chiarito fin da subito il diverso ruolo dei due elementi caratterizzanti il problema, necessario per riconoscere la stessa struttura in tutte e tre le fasi del progetto: vi sono elementi da cui si parte e si arriva, ed elementi che li "congiungono" e che quindi vengono attraversati.

2. L'attenzione al linguaggio, usato per rimarcare gli elementi del problema nelle sue varie forme, è stato di supporto alla *traduzione* della mappa di Königsberg nel relativo grafo: infatti, ha reso spesso spontanea l'identificazione delle regioni-vertici e dei ponti-archi. Ovviamente, non sono mancati episodi in cui la prima reazione è stata quella di identificare i ponti con i vertici: proprio in questi casi, è stato sufficiente rimarcare quali fossero gli elementi da "attraversare" per far convergere le classi verso una rapida e corretta traduzione.

3. Anche l'introduzione della rappresentazione simbolica dei cammini rappresenta un aspetto linguistico a cui si è prestata attenzione. Vi sono state classi in cui tale rappresentazione è emersa spontaneamente in uno o due partecipanti fin dalle prime due mappe, soprattutto nel caso delle classi del LM. Di norma però, la maggior parte dei presenti è sempre rimasta ancorata al linguaggio geometrico durante lo svolgimento del problema delle mappe.

## 5.3 Coreo-graph ed embodiment

1. Siamo stati consapevoli fin dall'inizio che proporre un'attività legata alla danza avrebbe comportato dei rischi: il rischio intrinseco che, riconoscendo le soluzioni direttamente sullo schema, ragazzi e ragazze non sarebbero stati invogliati a eseguire le coreografie; e il rischio che la danza avrebbe potuto creare in alcuni/e una situazione di disagio. La prima difficoltà è stata superata proprio dall'aver previsto in fase di progettazione la separazione tra la fase della ricerca delle soluzioni (e della loro rappresentazione simbolica) da quella motoria, invitandoli ad assimilare anzitutto gli schemi motori e successivamente il coordinamento con la base musicale. Rispetto alla seconda potenziale difficoltà, abbiamo sempre invitato a cambiare ballerino/a a ogni nuova coreografia Euleriana da riprodurre, lasciando ampia libertà nella personalizzazione degli stili dei passi. Nonostante alcune iniziali perplessità, abbiamo riscontrato un buon livello di coinvolgimento e abbiamo spesso osservato una collaborazione tra i membri di una stessa crew, in cui qualcuno ha dedicato maggior impegno nella ricerca delle coreografie sugli schemi, e qualcun altro nel verificare le soluzioni proposte ballando le coreografie. In particolare, abbiamo riscontrato una maggior reticenza nelle ragazze della SSPG rispetto ai ragazzi, che (contrariamente alle nostre attese indotte dai pregiudizi di genere sulla danza) si sono spesso mostrati più disinvolti e partecipi.

2. Il fatto di percorrere fisicamente i grafi (qui mascherati da schemi di danza) ha rafforzato l'osservazione di un'importante proprietà delle soluzioni: quando esistono coreografie Euleriane partendo da ogni vertice, queste sono sempre dei circuiti, mentre se esistono partendo solo da due precise posizioni dello schema, allora tali posizioni sono sempre l'inizio e la fine della coreografia Euleriana. Riteniamo che questo rinforzo sia stato prezioso per associare la differenza tra circuiti-cammini al differente ruolo svolto dalle posizioni di partenza. Tale associazione, ripresa successivamente nell'attività sui grafi, è stata l'argomento che abbiamo sviluppato per guidare le classi verso la corretta formulazione dell'enunciato.

3. Aver anticipato l'attività sui cammini Euleriani in un grafo, inserendola in un contesto motorio stimolante e coinvolgente, ha fornito a tutti gli studenti una base di esperienze a cui attingere per la parte concettualmente più faticosa del laboratorio, ovvero la formulazione del Teorema di Eulero. Ci siamo accorti che per una parte degli studenti (tra cui quelli con maggiori difficoltà a mantenere la concentrazione) l'attività su carta e penna proposta nella terza fase è risultata poco interessante o faticosa: per questi studenti, l'esperienza di dover "passare" fisicamente da una posizione a una differente con un passo è stato un riferimento concreto a cui rifarsi durante la terza fase. Tale esperienza è stata richiamata durante la spiegazione dell'argomento *entry-exit*, per riflettere sul valore del grado nel caso dei vertici di passaggio (ogni abbandono della posizione è compensato da un ritorno, e viceversa), rispetto ai vertici di inizio e fine della coreografia.

4. Rispetto al problema della guida turistica proposto in [Gaio & al. (2020)], la cui principale difficoltà consiste nella ricerca e svolgimento di un percorso Euleriano, abbiamo qui introdotto la a difficoltà aggiuntiva di dover sincronizzare la coreografia Euleriana con il ritmo della base musicale: questo comporta dover coordinare la percezione dello spazio e del tempo [Gerofsky (2013)], rispettando una schema geometrico preassegnato dato dalla coreografia. Consapevoli che ciò avrebbe potuto costituire un ostacolo all'attività di *embodiment*, e ispirandoci ad alcune metodologie dell'insegnamento della danza, abbiamo separato l'esecuzione delle sole coreografie Euleriane, come momento di verifica dell'analisi geometrica e memorizzazione dei soli schemi motori (e degli stili dei passi), dalla sincronizzazione con la base musicale, avvenuta in un secondo momento, una volta che gli schemi sono stati acquisiti.

### 5.4 Formulazione dell'enunciato: dal globale al locale e ritorno

1. La terza fase sintetizza tutto il percorso svolto nel progetto: dal riconoscimento di proprietà globali (la ricerca di cammini Euleriani sul grafo obbliga a visualizzare tutto il grafo), alla scoperta di aspetti locali rilevanti per lo studio globale (la parità del grado dei vertici), fino alla formulazione di un enunciato (il Teorema di Eulero). Di questi aspetti, sicuramente la formulazione dell'enunciato è la parte che ha spesso creato maggior difficoltà; infatti, avere a disposizione un numero molto ridotto di esempi a cui riferirsi, nonostante le tre fasi abbiano ripetuto sostanzialmente lo stesso schema, può condurre verso congetture non completamente corrette. Focalizzarsi su un numero limitato di esempi specifici può essere fuorviante, e comunque insufficiente a immaginarsi il caso generale di riferimento [Kapur (2010), Kapur (2011)]. Per ricondurre la narrazione verso la corretta congettura, eliminando di volta in volta falsi candidati, ci siamo quindi serviti di esempi aggiuntivi, come quelli in **Figura 7**.

2. Di fronte alla difficoltà di cogliere il ruolo discriminante del grado di un vertice, con riferimento al grafo G2 dell'allegato 1, abbiamo disegnato il cammino Euleriano come un vero e proprio tragitto, in analogia agli enigmi che chiedono di disegnare una casetta senza staccare la penna dal foglio (e richiamando l'esperienza *embodied*, nel caso del SSPG). Nel farlo, abbiamo insistito sui termini "uscire da" ed "entrare in" un vertice: così, da tutti i vertici di passaggio si entra e si esce lo stesso numero di volte, mentre dal vertice di partenza si esce una volta in più di quante non si rientri e nel vertice di arrivo si entra una volta in più di quante non se ne esca. L'uso di tali termini e le osservazioni sulla differenza tra vertici di passaggio e di partenza/arrivo, ha reso evidente il ruolo svolto dai gradi dei vertici rispetto ad altre proprietà locali/globali del grafo (parità/disparità di vertici o archi per esempio). In questo cruciale passaggio, che porta a riconoscere il ruolo del grado del vertice nell'esistenza dei cammini Euleriani, abbiamo riscontrato una sostanziale differenza tra SSPG e LM. Nelle classi della SSPG, il contenuto dell'enunciato del Teorema di Eulero ci è parso esser stato intuito solo da pochi, che infatti hanno poi guidato le classi nell'applicazione dell'enunciato al problema delle mappe. In queste classi, dopo aver ragionato sul grafo G2 come sopra descritto, ci siamo limitati a verificare la validità dell'enunciato, senza insistere troppo sull'argomento logico-dimostrativo. Diversamente, le classi del LM hanno spontaneamente cercato di estendere l'argomento usato sul grafo G2 agli altri tipi di grafo, con un minimo sforzo da parte nostra nel guidarli a intravedere uno schema di (parziale) dimostrazione dell'enunciato.

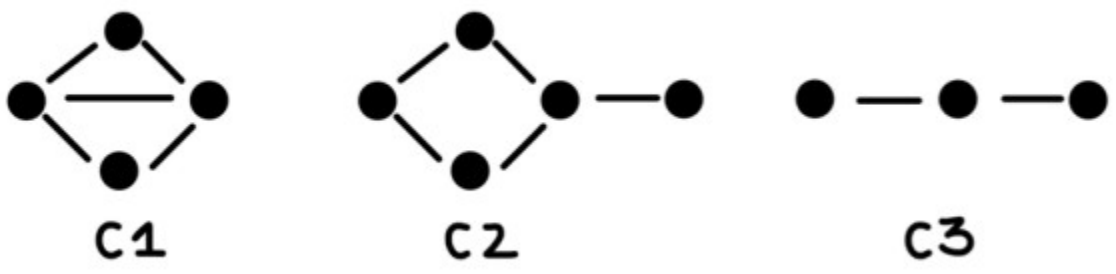

**Figura 7.** Controesempi: G2 e G3 hanno entrambi 5 vertici, mentre C1 ne ha solamente 4. In G2, G3 e C1 compaiono dei triangoli, mentre in C2 no. In tutti i grafi precedenti c'è un ciclo, in C3 no. Inoltre, C3 è l'unico tra i grafi presentati che ha un numero di vertici minore del numero di archi.

### 5.5 Problema diretto e problema inverso

In questo laboratorio abbiamo insistito molto sul far lavorare le classi anche sul problema inverso dell'esistenza dei cammini Euleriani, proponendolo in almeno due situazioni:

- la richiesta di come modificare gli schemi di danza per avere assegnate coreografie Euleriane (inizio della terza fase);
- la rimozione e ricostruzione di un ponte nell'attività conclusiva sul Teorema di Eulero (conclusione della terza fase).

In entrambi i casi, per quanto la prima richiesta sia stata proposta *prima* dell'enunciato del Teorema e la seconda *dopo* l'enunciato, l'obiettivo è sempre stato di cercare di indirizzare l'attenzione verso il grado del vertice, come proprietà locale discriminante. Lungo tutta la sperimentazione, entrambe le versioni del problema inverso sono sempre state accolte dalle classi come *sfidanti*, e hanno rappresentato un importante momento di sviluppo e comprensione dell'attività, anche quando la strategia risolutiva si è rivelata discosta dalle nostre aspettative: infatti, nello svolgimento dell'ultima attività sull'abbatimento e ricostruzione di un solo ponte, la maggior parte dei partecipanti (sia per la SSPG che per il LM) ha lavorato direttamente sulla mappa, piuttosto che con le operazioni di taglia/incolla di un arco sul corrispondente grafo, mostrando un'evidente preferenza verso una rappresentazione meno astratta del problema.

## 6 Conclusioni e sviluppi successivi

La struttura "circolare" del progetto, nella sua sequenza di varie fasi laboratoriali che conducono alla teoria dei grafi e alla risoluzione del problema dei ponti di Königsberg, ha confermato la sua efficacia come attività di *problem-solving*: essa ha permesso di coinvolgere le classi dalle prime battute, mantenendo viva la loro curiosità fino agli sviluppi finali e offrendo loro una prospettiva storico-sperimentale della disciplina a cui non sono abituati. Più in generale, questa esperienza ha rafforzato quanto già sostenuto nella letteratura circa le potenzialità didattiche della teoria dei grafi per sviluppare nella SSPG il pensiero critico e l'attitudine verso il *problem-solving* [Mammana & Milone (2010), Ferrarello & Mammana (2017), Robinson (2006), Santoso (2018)].

Le nostre modifiche ed integrazioni al progetto contenuto in [Gaio & al. (2020)], volte a mettere maggiormente in luce le caratteristiche dei cammini Euleriani e il ruolo del grado dei vertici, si sono mostrate utili al fine di far emergere l'enunciato del Teorema di Eulero come una naturale sintesi delle esperienze laboratoriali, capace a posteriori di risolvere situazioni ben più complesse di quelle assegnate. Come emerso dai questionari somministrati nel corso di questo triennio, questo progetto ha rinforzato il valore della matematica come linguaggio che svela la struttura di un problema ed è in grado di fornire strumenti tanto inattesi quanto efficaci per le strategie di *problem solving*.

Tali modifiche si sono mostrate particolarmente utili nel caso della sperimentazione con le classi del biennio del LM. In questo caso, il progetto è stato svolto direttamente presso il Dipartimento di Matematica[3] nell'arco di mezza giornata, senza la possibilità di includere alcuna attività di tipo motorio. A fronte di uno sviluppo del progetto ridotto e condensato in poche ore consecutive, il risultato è stato comunque molto positivo: le classi del LM hanno spesso mostrato un ottimo livello di coinvolgimento, adattandosi con facilità agli aspetti più formali delle attività (rappresentazione simbolica, teoria dei grafi), talvolta anticipando in maniera intuitiva alcuni aspetti geometrici della struttura del problema (chiusura o apertura dei cammini, ruolo del grado di un vertice).

[3] Dipartimento F.Enriques, Univ. Statale di Milano.

L'attività motoria *Coreo-Graph*, esclusa dalla sperimentazione nel LM, ha invece evidenziato i suoi effettivi vantaggi con le classi della SSPG. Essa:

- ha permesso di anticipare in un contesto motorio-esperienziale gli oggetti e i quesiti della teoria dei grafi, proposti successivamente in classe, facilitandone lo studio più astratto. Tale attività ha diluito le difficoltà della terza fase, che altrimenti sarebbe stata troppo densa e faticosa, e ha offerto un costante riferimento tangibile ai concetti della teoria dei grafi che hanno condotto all'enunciato del Teorema di Eulero.
- Ha fatto sperimentare alle classi un'attività matematica in un contesto differente dalla solita didattica frontale in aula, ampliando l'interesse e l'efficacia didattica all'interno di classi di norma molto eterogenee. Infatti, è capitato che chi si sia trovato adeguatamente stimolato dalle attività dei ponti e dei grafi, sia stato "disorientato" dalla natura motoria di questo nuovo compito. E viceversa, che qualcuno con difficoltà di concentrazione e di apprendimento, o poca motivazione al lavoro in classe, abbia trovato maggiori stimoli in questo tipo di attività.

Questa sperimentazione ci suggerisce che l'attività motoria *Coreo-Graph* sia stata didatticamente efficace e significativa per le classi della SSPG, andando incontro alla richiesta di un'esperienza capace di creare idee matematiche attraverso la pratica laboratoriale: essa ha infatti permesso di ancorare un concetto astratto, quello di *cammino Euleriano su un grafo connesso*, a una comune esperienza motoria nello spazio. In tal senso, crediamo che la dimensione del mesospazio che ha caratterizzato la nostra proposta (o altre equivalenti che si possono immaginare e sperimentare) sia stata più efficace rispetto alla scelta originale di far lavorare le classi all'interno di un percorso urbano, dove l'osservazione della struttura del cammino e delle sue proprietà è meno diretta.

Di contro, per le classi del biennio del LM l'attività iniziale sulle mappe e i ponti sembra esser stata sufficientemente propedeutica all'attività teorica sui grafi, probabilmente grazie alla maggior attitudine verso la disciplina e a una capacità di astrazione già più sviluppata in questa fascia d'età. In questo non vediamo necessariamente un ridimensionamento dell'utilità didattica di questa o di altre attività *embodied* nelle classi della Scuola Secondaria di Secondo Grado.

Infatti, si sarebbe potuto disegnare il progetto diversamente, anticipando l'attività *Coreo-Graph* come preparatoria alla parte teorica, e spostando l'applicazione al problema dei ponti di Eulero in coda al progetto. Una prima indicazione del possibile valore di questa variante viene dalle prime due edizioni dell'evento "Un anno di Licei Matematici" [4], svoltesi a conclusione dell'anno scolastico del LM (Maggio 2023 e 2024). In tale occasione, abbiamo proposto ad alcuni gruppi del LM direttamente il problema delle coreografie Euleriane, a cui abbiamo fatto seguire la formulazione del Teorema di Eulero, usando gli schemi di danza come unici esempi rappresentativi dei casi presenti nell'enunciato.

Rimane dunque aperta la domanda circa la significatività di un *laboratorio embodied* per le classi della Scuola Secondaria di Secondo Grado, sul quale cercheremo di concentrare i nostri ed altrui sviluppi futuri di questo progetto.

Oltre a portare avanti la sperimentazione con queste fasce d'età all'interno degli Istituti Scolastici del nostro territorio, valutandone possibili estensioni nell'area della Computer Science, per il futuro vorremmo consolidare l'attività *Coreo-Graph* come laboratorio indipendente. Lo stiamo infatti già proponendo come strumento per diffondere delle attività di *problem solving* anche al di fuori del contesto strettamente scolastico, inserendolo nella programmazione di alcuni Centri Estivi di danza/multisport.

## Ringraziamenti

Ringraziamo Laura Branchetti, per averci stimolato ed incoraggiato nella sperimentazione e stesura di questo progetto. Ringraziamo Alessia Maggioni, Maria Luoni, Camilla Pugliese, Giorgia Mazzilli e Beatrice Sacchetti per l'aiuto durante gli incontri con le classi dei Licei Matematici che si sono susseguite in questi anni. Ringraziamo le insegnanti e le classi dell'Istituto Comprensivo "A. Manzoni" di Ornago e Burago per aver accolto questa proposta laboratoriale, e le classi del Liceo Matematico che hanno chiesto e partecipato con entusiasmo al laboratorio di teoria dei grafi nel periodo 2022-2024. T.P. ringrazia infine Alexandar Tasovac detto Sasko, della scuola di danza Danza Passion di Burago di Molgora (MB), per il supporto concreto ed i piccoli, ma preziosi, consigli creativi sull'attività *Coreo-graph*.

[4] https://liceomatematicomilano.unimi.it/

## Riferimenti bibliografici

# Allegato 1

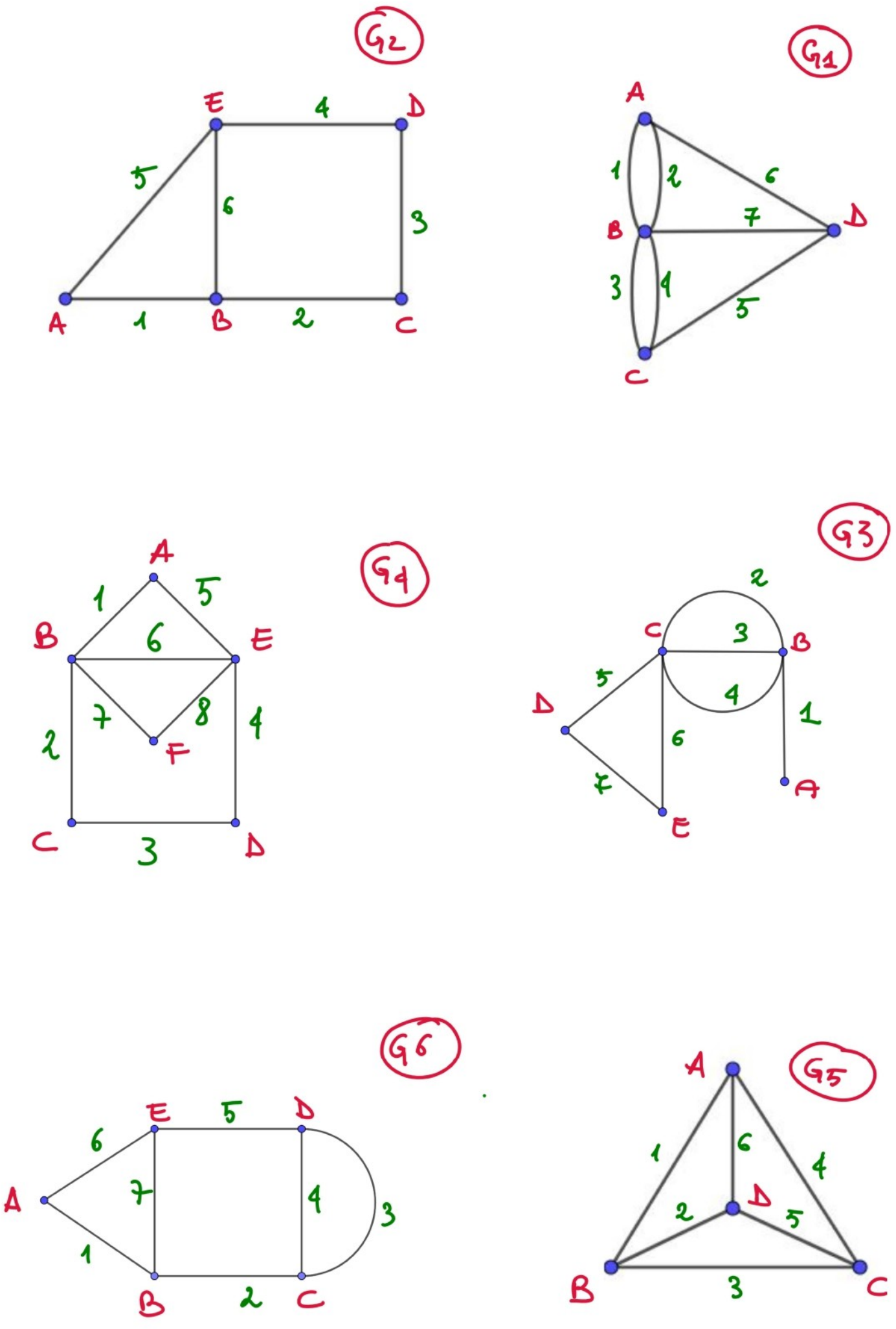

**Figura 8.** Scheda per la ricerca di cammini Euleriani su grafi usata nella terza fase